\pdfoutput=1

\documentclass[pdflatex,sn-mathphys-num]{sn-jnl}

\usepackage{hyperref}
\hypersetup{
    unicode      = false,
    pdftoolbar   = true,
    pdfmenubar   = true,
    pdffitwindow = true,
    pdftitle     = {Parallel MADS},
    pdfauthor    = {Le~Digabel, Lesage-Landry, Mendoza, Tribes},
    pdfsubject   = {Parallel MADS},
    pdfnewwindow = true,
    pdfkeywords  = {Parallelism, MADS},
    colorlinks   = true,
    linkcolor    = blue,
    citecolor    = blue,
    filecolor    = black,
    urlcolor     = blue,
    breaklinks   = true
}

\usepackage[english]{babel}
\usepackage{anyfontsize} %
\usepackage{multirow}
\usepackage{amsmath}
\usepackage{amssymb}
\usepackage{amsfonts}
\usepackage{amsthm}
\usepackage{mathrsfs}
\usepackage[title]{appendix}
\usepackage{xcolor}
\usepackage{textcomp}
\usepackage{manyfoot}
\usepackage{booktabs}
\usepackage[vlined, ruled, linesnumbered]{algorithm2e}
\usepackage{geometry}
\usepackage{xspace}
\usepackage{tikz}
\usetikzlibrary{fit, positioning, calc, decorations.pathreplacing}
\usepackage{graphicx}
\usepackage{float}
\usepackage{dsfont}

\usepackage{pgfplots} %
\pgfplotsset{
  compat=1.18,   
  tick label style = {font=\rmfamily},
  every axis label = {font=\rmfamily},
  legend style = {font=\rmfamily},
  label style = {font=\rmfamily}
}

\newcommand*{\tikzmk}[1]{\tikz[remember picture,overlay,] \node (#1) {};\ignorespaces}
\newcommand{\boxit}[1]{\tikz[remember picture,overlay]{\node[yshift=3pt,fill=#1,opacity=.25,fit={(A)($(B)+(.90\linewidth,.8\baselineskip)$)}] {};}\ignorespaces}

\usepackage{soul} %
\usepackage{colortbl} %

\definecolor{Red}{rgb}{1,0,0}
\definecolor{Gray}{rgb}{0.2,0.2,0.2}
\definecolor{Maroon}{rgb}{0.6,0.05,0.03}
\definecolor{Blue}{rgb}{0,0.7,0.9}
\definecolor{Green}{rgb}{0,.7,0}

\def\R{{\mathbb{R}}}
\newcommand{\xx}{\mathbf{x}}
\newcommand{\X}{\mathcal{X}}
\newcommand{\PP}{\mathds{P}}

\newcommand{\bbo}{{BBO}\xspace}
\newcommand{\dfo}{{DFO}\xspace}
\newcommand{\ds}{{DS}\xspace}
\newcommand{\mads}{{\sf MADS}\xspace}
\newcommand{\nomad}{{\tt NOMAD}\xspace}

\usepackage[nameinlink,capitalise,noabbrev]{cleveref}
\usepackage{silence}
\usepackage{caption}
\usepackage{subcaption}

\graphicspath{{figures/}}
\usepackage{comment}

\makeatletter
\def\doiurl#1{%
  \url{https://doi.org/#1}
  \edef\sn@thisdoi{https://doi.org/#1}%
  \sn@lookahead}
\def\sn@lookahead{\@ifnextchar.{\sn@eatdot}{\@ifnextchar\burl{\sn@checkburl}{}}}
\def\sn@eatdot.{\@ifnextchar\burl{\sn@checkburl}{ . }}
\def\sn@checkburl\burl#1{%
  \edef\sn@tmp{#1}%
  \ifx\sn@tmp\sn@thisdoi\else\ . \sn@burl{#1}\fi}
\def\sn@burl#1{\textsf{#1}}
\def\burl#1{\sn@burl{#1}}
\makeatother

\begin{document}

\title[Parallel versions of the MADS algorithm]{Parallel versions of the mesh adaptive direct search algorithm}

\author[2,3]{\fnm{S\'ebastien} \sur{Le~Digabel}}\email{\tt sebastien.le-digabel@polymtl.ca}
\author[1,3]{\fnm{Antoine} \sur{Lesage-Landry}}\email{\tt antoine.lesage-landry@polymtl.ca}
\author[1,3]{\fnm{Samuel} \sur{Mendoza}}\email{\tt samuel.mendoza@polymtl.ca}
\author[2,3]{\fnm{Christophe} \sur{Tribes}}\email{\tt christophe.tribes@polymtl.ca}

\affil[1]{\orgdiv{Department of Electrical Engineering}, \orgname{Polytechnique Montr\'eal}, \orgaddress{\street{2500 Chemin de Polytechnique}, \city{Montr\'eal}, \postcode{H3T~1J4}, \state{Qu\'ebec}, \country{Canada}}}

\affil[2]{\orgdiv{Department of Mathematics and Industrial Engineering}, \orgname{Polytechnique Montr\'eal}, \orgaddress{\street{2500 Chemin de Polytechnique}, \city{Montr\'eal}, \postcode{H3T~1J4}, \state{Qu\'ebec}, \country{Canada}}}

\affil[3]{\orgname{Group for Research in Decision Analysis (GERAD)}, \orgaddress{\street{3000 Chemin de la C\^ote-Sainte-Catherine}, \city{Montr\'eal}, \postcode{H3T~2A7}, \state{Qu\'ebec}, \country{Canada}}}

\abstract{This work surveys the different parallel variants of the mesh adaptive direct search (\mads) algorithm for constrained blackbox optimization. These problems can inherently imply high computational costs due to the possible large number of variables and multi-modality of the search space.
In addition, the potential time-intensive nature and time heterogeneity of the blackboxes defining the problem prompts the need for efficient implementations. Parallelism emerges as an actionable solution to mitigate computation time, as modern computer systems rely on multi-core architecture. The reviewed methods employ diverse levels of parallelism and distinct parallel strategies to effectively tackle each aspect outlined above. The manuscript details the practical implementations, provides computational results, and offers insights into the advantages and limitations of each \mads parallel method.}


\maketitle

\section{Introduction}

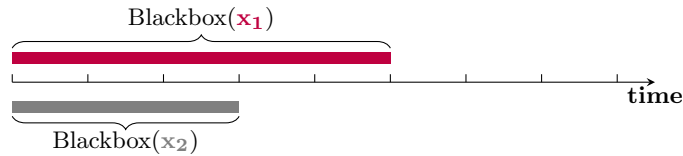
\begin{figure}[ht!]
\centering
\begin{tikzpicture}[%
every node/.style={
font=\scriptsize,
text height=1ex,
text depth=.25ex,
},
]
\draw[-stealth] (0,0) -- (8.5,0);

\foreach \x in {0,1,...,8}{
\draw (\x cm,3pt) -- (\x cm,0pt);
}
\node[anchor=north] at (8.5,0) {\bf\textrm{time}};
\fill[purple] (0,0.25) rectangle (5,0.4);
\fill[gray] (0,-0.25) rectangle (3,-0.4);
\draw[decorate,decoration={brace,amplitude=5pt}] (0,0.45) -- (5,0.45)
node[anchor=south,midway,above=4pt] {\textrm{Blackbox(\textcolor{purple}{$\mathbf{x_1}$})}};

\draw[decorate,decoration={brace,amplitude=5pt}] (3,-0.45) -- (0,-0.45)
node[anchor=north,midway,below=4pt] {\textrm{Blackbox(\textcolor{gray}{$\mathbf{x_2}$})}};

\end{tikzpicture}
\caption{Illustration of a time-heterogeneous blackbox.}
\label{fig-heter-bb}
\end{figure}

This work considers optimization problems of the form
\begin{equation}\tag{$\PP$}
\label{pb-def}
\min_{\xx \in \Omega} f(\xx),
\end{equation}
where $\Omega = \left\{\left. \xx \in \X ~\right|~ c_{j}(\xx) \leq 0, j \in \mathcal{J}=\{1, 2, \ldots, m\} \right\} \subseteq \R^n$, $m, n \in \mathbb{N}$, is the feasible set and $\X$ is a subset of $\R^n$, typically defined by bound constraints.
The functions $f: \X \rightarrow \R \cup \{ \infty \}$ and $c_{j}: \X \rightarrow \R \cup \{ \infty \}$, $j \in \mathcal{J}$,
are typically evaluated through a computer simulation with no available derivatives,
and is considered as a blackbox.
We consider that evaluating such a blackbox is time-consuming and that only a limited budget of evaluations is available. We further assume that the blackbox may be subject to heterogeneity in the time required to evaluate a point, as illustrated in \cref{fig-heter-bb}.

Problem~\eqref{pb-def} is a blackbox optimization (\bbo) problem, for which derivative-free optimization (\dfo) algorithms are considered. These techniques are described in~\cite{AuHa2017,CoScVibook},
and some practical applications are illustrated in~\cite{AlAuGhKoLed2020}.

With \dfo methods, a rule of thumb is to allow a budget of the order of 1,000$n$ evaluations
to obtain satisfactory results.
This is not conceivable when blackboxes are time-consuming because this process would lead to impractical resolution times.
The situation worsens when dealing with large-scale problems, typically with more than 50 variables, as the evaluations budget escalates significantly.

Parallel computing appears as a solution because most of computational resources exploits multi-core architectures available through various levels, e.g.,  workstations, servers, and supercomputers.
In addition, with the end of Dennard scaling, parallelization has become much more important than before, as the alternative of increasing processor clock speed is no longer possible due to heat dissipation and power consumption issues.

An obvious way of using parallelism is to ``open'' the blackbox and to modify it to exploit the available parallel resources.
However, this is not always possible due to many reasons: The source code of the blackbox may be unavailable, too complicated, too old, etc.
Hence, in this work, we focus on the parallelization of \dfo algorithms, and in particular direct search (\ds) methods.

In an effort to take advantage from modern multi-core infrastructures, multiple \dfo methods have been extended to leverage parallel computations through their resolution process.
In particular, \dfo methods can benefit from parallelism to reduce the time resolution, improve solution quality, and increase the size of \bbo problems that can be practically solved. 
Specifically, \ds methods, a branch of \dfo methods, provides a framework for the effective use of parallel computing.

Throughout this work, we consider the mesh adaptive direct search (\mads)~\cite{AuDe2006} algorithm and its {\sf C++} implementation, \nomad version~3~\cite{Le09b}, and we review its different parallel versions.
We measure the performance of these parallel algorithms by introducing metrics for benchmarking in a context of parallel computing, and we give recommendations to guide users on what kind of parallelism should be used depending on the different aspects of the problem at hand.
This study is conducted with \mads, but remains applicable to other \ds methods.

This work is structured as follows: \cref{literature} reviews the literature and provides an overview of parallel methods for \bbo,
\cref{sec-mads} describes \mads,
\cref{parallel_nomad} presents its parallel versions,
and \cref{sec-tests} reports computational experiments with recommendations. Finally, concluding remarks are provided in \cref{conclusion}.

\section{Parallelism in \bbo}
\label{literature}

A guiding principle for the efficient use of parallelization to solve~\eqref{pb-def} is to consider that the time spent to evaluate the objective or constraint functions is much higher than the time required by an optimization algorithm to generate trial points. Hence, the parallelization that is the most likely to be beneficial, i.e., which yields a large reduction in overall execution time,
consists of performing several concurrent evaluations of the blackbox and/or to perform an internal parallelization of the blackbox~\cite{Olsson2014}. As mentioned above, we consider that the latter is not available in this work.
Parallelization strategies, i.e., functional parallelism, domain decomposition, multi-search parallelism, and the hybridization of mutually compatible aforementioned strategies, are described in details in~\cite{Sch2020}. It also provides a literature review on integrative framework for parallel computational optimization in operations research. While these strategies are formulated in an operations research context, the scope of applicability can be straightforwardly extended to other optimization algorithm classes, such as \ds methods, as described in~\cite{DeWu}.

Parallelism in \bbo and specifically for \ds methods, is a longstanding topic of interest. The first \ds methods to employ parallelism with proofs of convergence are based on the generalized pattern search ({\sf GPS})~\cite{Torc97a} and the generating set search ({\sf GSS})~\cite{KoLeTo03a} algorithms.
The resulting method is the asynchronous parallel pattern search~\cite{GrKo2010,GrKoLe2008,HoKoTo01a,Ko05,KoTo04a}. Implementation examples can be found in the {\tt APPSPACK}~\cite{GrKo06} and {\tt HOPSPACK}~\cite{Hops20-Sandia} solvers.
Reference~\cite{LiTr2017} combines instances of \mads with different line search methods which are run in parallel to accelerate the resolution process. A specific parallel search for \mads is later introduced in the technical report~\cite{TaAlKo2020}.
Other \ds methods include the Nelder-Mead~\cite{NeMe65a} whose first parallel version is described in~\cite{DeTo91a}, and revisited more recently in~\cite{OzWaOn2019}. The {\sf DIRECT} algorithm~\cite{JoPeSt93a} is hybridized with {\sf GSS} in~\cite{griffin2010asynchronous} and its parallelized versions are analyzed in~\cite{He2009,He2007,He2009a,He2009b}.
Other \dfo methods include model-based algorithms utilizing a trust region and a model such as polynomial or radial basis functions. In this case, parallelization is not as inherent as for \ds methods. Nonetheless, examples of parallel \dfo can be found in~\cite{FVBerghen_2004,FVBerghen_HBersini_2005,HoMe02a,ReSh2007,CASRGR03,Xia2021}.
Parallel surrogate-assisted methods are also described in~\cite{briffoteaux23,gagagaroco20,Haftka16}.
Bayesian optimization methods based on the efficient global optimization framework~\cite{JoScWe1998} are parallelized, for example, in~\cite{TRAN2019827,Zhan2017,Zhan2020}. 
There exists also multiple other parallel heuristics such as~\cite{VAZQUEZ2016599}. More generally, the survey~\cite{Alba13} reviews different kinds of parallel heuristics.

Finally, parallel \ds methods are employed in specific applications such as hyperparameters tuning~\cite{AuDaOr13a}, local optimum enumeration~\cite{LW16}, multiobjective optimization~\cite{Tavares2023}, and small (or moderate) dimensional engineering design problems driven by computationally expensive simulations, where large-batch parallel sampling can outperform more elaborate techniques~\cite{PANG2023159544,CAMPOS2025109422}.

\section[MADS: Mesh Adaptive Direct Search]{\mads: Mesh Adaptive Direct Search}
\label{sec-mads}

\mads is a direct search method introduced in~\cite{AuDe2006} and refined in several articles, including its last major evolution in~\cite{AuLeDTr2018}. The version given in~\cref{algo:mads} corresponds to the one provided in~\cite{AuHa2017}.

\begin{algorithm}[ht]
  \caption{\mads with Extreme Barrier}
  \label{algo:mads}
  \SetAlgoLined
  \SetKwBlock{Init}{Initialization}{}
  \Init{
    \hspace*{1.65mm}Set $\Delta^{0} > 0$ and $\xx^0 \in \Omega$. \\
    Start the $w$ worker processes for evaluation. \\
    Evaluate $\{f(\xx^0), c_1(\xx^0), c_2(\xx^0), \dots, c_m(\xx^0)\}$ and initialize cache $\mathcal{C}$.
  }
  \For{$k = 0, 1, \dots$}{
  \SetKwBlock{Mesh}{Mesh {\normalfont(construction of $\mathcal{M}^k$~\eqref{eq-mesh})}}{}
  \Mesh{  
  \hspace*{1.65mm}Set $\mathbf{D} \in \mathbb{R}^{n \times p}$
  a positive spanning matrix. \\
  Set $\delta^k = \min\{ \Delta^k, (\Delta^k)^2 \}$. \\
  }
  \SetKwBlock{Search}{Search {\normalfont(optional)}}{}
  \Search{
  \hspace*{1.65mm}Build the search set $\mathcal{S}^k\subset \mathcal{M}^k$. \\
  Evaluate $\{f(\xx), c_1(\xx), c_2(\xx), \dots, c_m(\xx)\}$ \text{ for all } $\xx \in \mathcal{S}^k$\\
  \hspace*{1.65mm}(opportunistically)\\
  \If{$f_{\Omega}(\mathbf{x}) < f_{\Omega}(\mathbf{x}^k)$ \normalfont{for some} $\mathbf{x} \in \mathcal{S}^k$ }{ set $\mathbf{x}^{k+1} \leftarrow \mathbf{x}$ , $\Delta^{k+1} \leftarrow \tau^{-1}\Delta^k$\\go to \textbf{Termination}.}
  }
  \SetKwBlock{CacheSearch}{Cache Search {\normalfont(optional and only for $k > 0$)}}{}
  \CacheSearch{
    \If{$f_{\Omega}(\xx^*) < f_{\Omega}(\xx^k) \; \mathrm{for} \; \mathrm{some} \; \xx^* \in \mathcal{C}$}{set $\xx^{k+1} \leftarrow \xx^*$ , $\Delta^{k+1} \leftarrow \tau\Delta^*$ \\ go to \textbf{Termination}.}}
  
  \SetKwBlock{Poll}{Poll}{}
  \Poll{
   \hspace*{1.65mm}Build the poll set
   $\mathcal{P}^k \subset \mathcal{M}^k$.\\
   \tikzmk{A}
   Perform evaluations and compute success\\
   \If{Success}{go to \textbf{Termination}.}
   \hspace*{1.65mm}Set $\mathbf{x}^{k+1} \leftarrow \mathbf{x}^k$, $\Delta^{k+1} \leftarrow \tau\Delta^k$ (no success).\\
   \tikzmk{B}\boxit{gray}
   }
   \SetKwBlock{term}{Termination}{}
 \term{
 \If{Stopping criteria}{
 \hspace*{1.65mm}End the algorithm. \\ End the $w$ workers.}
 }
 }

\end{algorithm}

Blackbox evaluations occur at two different instances: The search and the poll.
The points to be evaluated are gathered in two sets $\mathcal{S}^k \subset \R^n$ and $\mathcal{P}^k \subset \R^n$ and are evaluated opportunistically, i.e., evaluations are interrupted as soon as a new best point is found.
When evaluations are computationally costly, as in \bbo, this strategy is effective because it can reduce the number of blackbox evaluations which is usually subject to a budget constraint. In a sequential programming environment, opportunism occurs at an individual point level. Thus, a minimum of one evaluation can be computed in order to terminate the current iteration if a simple decrease condition is reached by an evaluation point. In contrast, in a parallel environment, evaluations are dispatched concurrently to several computing resources. The minimum number of evaluations that will be performed before the premature termination of an iteration is, therefore, equal to the number of computing resources being used in parallel. As a result, the evaluation budget is spent more quickly. This is why a computation time budget is a better-suited constraint in a parallel environment.

\mads relies on a special structure called the mesh:
at iteration $k$, it is a discretization of the space of variables
defined by
\begin{equation}
    \mathcal{M}^k = \left\{\left. \mathbf{x}^k + \delta^k \mathbf{D} \mathbf{y} ~\right|~ \mathbf{y} \in \mathbb{N}^p \right\}\subset \mathbb{R}^n,
\label{eq-mesh}
\end{equation}
where $\delta^k > 0$ is the mesh size parameter, $\mathbf{D} \in \R^{n \times p}$ is a positive spanning matrix, and $p \in \mathbb{N}$. 
\mads selects the candidate trial points on the mesh.

The notion of locality is parameterized by $\Delta^k > 0$, the frame size parameter, which defines the maximal distance between the incumbent solution $\mathbf{x}^k$ and the poll points  $\mathbf{x} \in \mathcal{P}^k$
with
$$\|\mathbf{x}-\mathbf{x}^k\|_{\infty} \leq b \Delta^k, $$
where $b=\max_{d' \in \mathbb{D}}\|d'\|_{\infty}$ and $\mathbb{D}$ is the set of columns of matrix $\mathbf{D}$.

At every iteration $k$, $\delta^k$ and $\Delta^k$ are resized via $\tau \in (0,1) \cap \mathbb{Q}$, the mesh size adjustment parameter, such that $0 < \delta^k \leq \Delta^k$ for all $k$. By doing so, a variety of directions that grow densely in the unit sphere is guaranteed to be generated.

All these conditions lead to the poll step being rigidly defined, but they ensure the convergence of the method. 
As detailed in~\cite{AuHa2017}, upon some assumptions about the problem,
global convergence of \mads to a local optimum is guaranteed.

The optional search step is more flexible as it allows any strategy to be used
for the construction of the sets $\mathcal{S}^k$. It can be generic or specific to a problem,
and only needs to generate points on the current mesh.
Examples of search strategies can be found in~\cite{AuBeLe08,AuTr2018,CoLed2011}.

\mads is often run with a cache to provide a history of the evaluations that have been carried out. The main purpose of the cache is to ascertain whether the point has been previously evaluated, thus preventing redundant evaluations.
Optionally, the cache can be used during a special search step only when trial points are added to the cache asynchronously (see \cref{sec-pmadsa}). If a point $\xx^*$ in the cache is better than $\xx^k$, it would be used as a new incumbent to start an iteration, and the corresponding frame size $\Delta^*$ is adopted. This step is disabled in a sequential environment.

In practice, the stopping criteria indicated in the algorithm can be multiple. Typically, they include a maximum number of blackbox evaluations, a time budget,  and/or a minimal mesh or frame size.

\mads provides two distinct approaches to deal with inequality constraints.
First, a straightforward extension with extreme barrier (EB), which simply defines
$f_{\Omega}(\xx)=f(\xx)$ if $\xx \in \Omega$, and 
$f_{\Omega}(\xx)=\infty$ otherwise.
\cref{algo:mads} is based on the EB.
The second approach is referred to as the progressive barrier (PB) and is defined in~\cite{AuDe09a}.
It considers a function $h(\xx)$ that measures the constraint violations at $\xx$
and ranks the iterates in a graph (called a filter) where compromises in terms of $f$ (objective) and $h$ (constraints) are defined.

Finally, \mads is available in the open-source solver \nomad~\cite{Le09b} at
\href{https://www.gerad.ca/en/software/nomad/}{\tt https://www.gerad.ca/en/software/nomad/}. This solver will be used for our numerical comparisons.

\section[Parallel versions of MADS]{Parallel versions of \mads}
\label{parallel_nomad}

This section describes the four parallel methods of \mads, namely 
p\mads-S for ``Parallel \mads-Synchronous'' (\cref{sec-pmadss}),
p\mads-A for ``Parallel \mads-Asynchronous'' (\cref{sec-pmadsa}),
\textsc{COOP}-\mads for ``Cooperative \mads'' (\cref{sec-coopmads}),
and
\textsc{PSD}-\mads for ``Parallel Space Decomposition \mads'' (\cref{sec-psdmads}).

These methods are defined such that the convergence properties of \mads also hold.
While \textsc{PSD}-\mads has been the subject of~\cite{AuDeLe07}, the other three methods have never been described or benchmarked in the literature.
These methods follow the master-worker paradigm typically implemented with
the \texttt{Message Passing Interface} (\texttt{MPI}) library~\cite{SnOtHuWaDo95a}.
Following the parallel strategies described in \cref{literature}, p\mads-A and p\mads-S are classified as functional parallelism, while \textsc{PSD}-\mads is classified as domain decomposition parallelism. Lastly, \textsc{COOP}-\mads is classified as multi-search parallelism.

\subsection{p\mads}

\begin{figure}[ht!]
  \centering
  \includegraphics[width=0.5\textwidth]{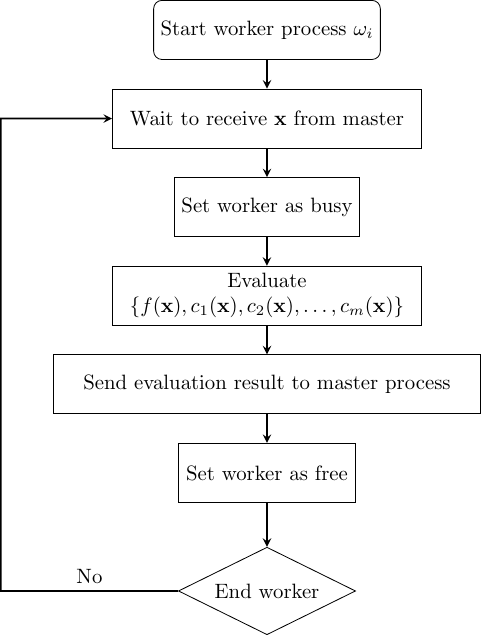}
  \caption{Flow chart of a single worker process for p\mads.}
  \label{fig:worker_pMads}
\end{figure}

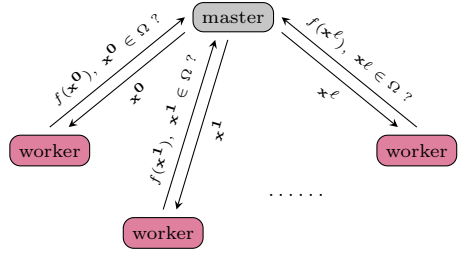
\begin{figure}[ht!]
     \centering
        \begin{tikzpicture}[every text node part/.style={font=\footnotesize}]
        
        \node[draw, rounded corners=4pt, fill=darkgray!30](master) {\textrm{master}};
        \node[draw,rounded corners=4pt, fill=purple!50](w1) at (180 + 36:30mm){\textrm{worker}};
        \draw[-stealth]([shift={(-1mm,0mm)}]master.south west)--([shift={(0mm,1mm)}]w1.40)node[midway,below, sloped]{ \tiny $\mathbf{x^0}$};
        \draw[-stealth]([shift={(0mm,1mm)}]w1.north)--([shift={(-1mm,0mm)}]master.180) node[midway,above, sloped]{ \tiny $f(\mathbf{x^0}), ~ \mathbf{x^0} \in \Omega~?$};

        \node[draw, rounded corners=4pt, fill=purple!50](w2) at (180 + 72:30mm){\textrm{worker}};
        \draw[-stealth]([shift={(0mm,-1mm)}]master.250)--([shift={(0mm,1mm)}]w2.50)node[midway,below, sloped]{ \tiny $\mathbf{x^1}$};
        \draw[-stealth]([shift={(0mm,1mm)}]w2.north)--([shift={(0mm,-1mm)}]master.220) node[midway,above, sloped]{ \tiny $f(\mathbf{x^1}), ~ \mathbf{x^1} \in \Omega~?$};

        \node[](w3) at (180 + 110:25mm){$\dots\dots$};

        \node [draw, rounded corners=4pt, fill=purple!50](w4) at (180 + 144:30mm){\textrm{worker}};
        \draw[-stealth]([shift={(0mm,1mm)}]w4.north)--([shift={(1mm,0mm)}]master.0) node[midway,above, sloped]{ \tiny $f(\mathbf{x^\ell}), ~ \mathbf{x\ell} \in \Omega~?$};
        \draw[-stealth]([shift={(1mm,0mm)}]master.south east)--([shift={(0mm,1mm)}]w4.140) node[midway,below, sloped]{ \tiny $\mathbf{x^\ell}$};
        \end{tikzpicture}
        \caption{Master-worker communication organization for p\mads.}
        \label{fig:pmads}
    \end{figure}

The first two parallel versions of \mads are both labelled p\mads and follow \cref{algo:mads}.
In these methods, the evaluation of trial points $\mathcal{P}^k$ is performed in parallel in a master-worker fashion (see \cref{fig:pmads}).
In doing so, the master is responsible for sending evaluation points to the workers, who are then tasked with evaluating them (see \cref{fig:worker_pMads}). The two versions of p\mads are
the synchronous variant (p\mads -S) and
the asynchronous variant (p\mads-A). They differ in the management of the evaluations by the $w$ workers (shaded area in \cref{algo:mads}). Starting and ending the workers are managed by the master process when the algorithm ends. 

\subsubsection{p\mads-S}
\label{sec-pmadss}

\begin{figure}[ht!]
  \centering
  \includegraphics[width=0.8\textwidth]{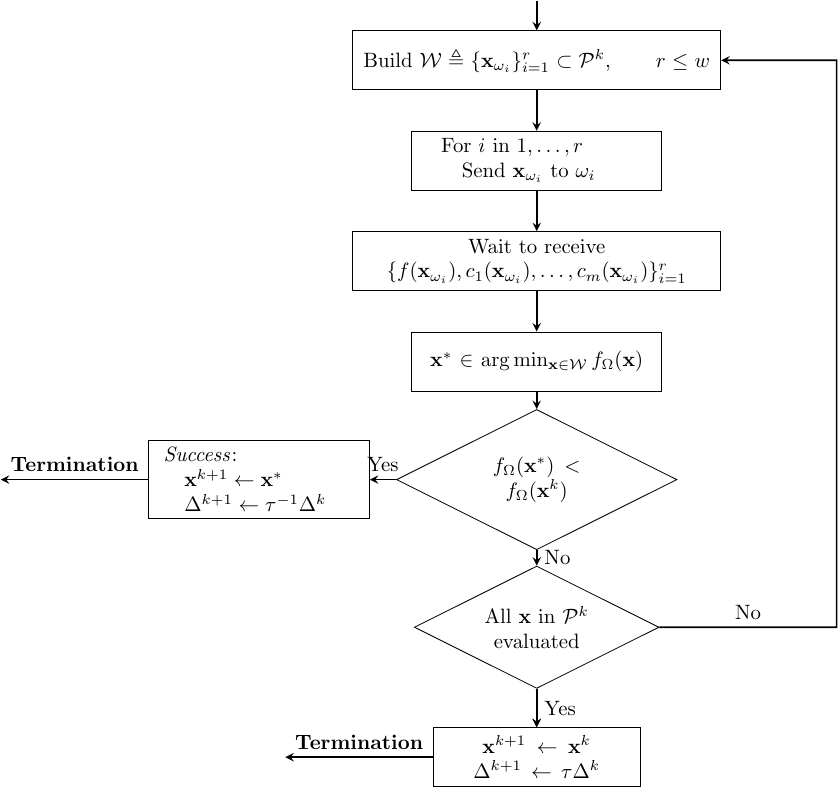}
  \caption{Master process flow chart of p\mads-S to manage parallel evaluations by the $w$ parallel workers allocated to \nomad, where $\omega_i$ represents their respective indices. Details the shaded area in Poll step of \cref{algo:mads}.}
  \label{fig:master_pMadsS}
\end{figure}

p\mads-S has a deterministic behaviour, but at the cost of introducing a synchronization barrier which restricts the master to end an iteration only when all evaluations in progress are terminated.  Consequently, this leads to idle computing resources at each iteration and hence to a loss in efficiency. The synchronization barrier eliminates the need for the cache search step.

An opportunistic block evaluation approach, using a block of size $w$, is achieved by executing the \textbf{Termination} step upon success, as shown in \cref{fig:master_pMadsS}. 

\subsubsection{p\mads-A}
\label{sec-pmadsa}

\begin{figure}[ht!]
  \centering
  \includegraphics[width=0.7\textwidth]{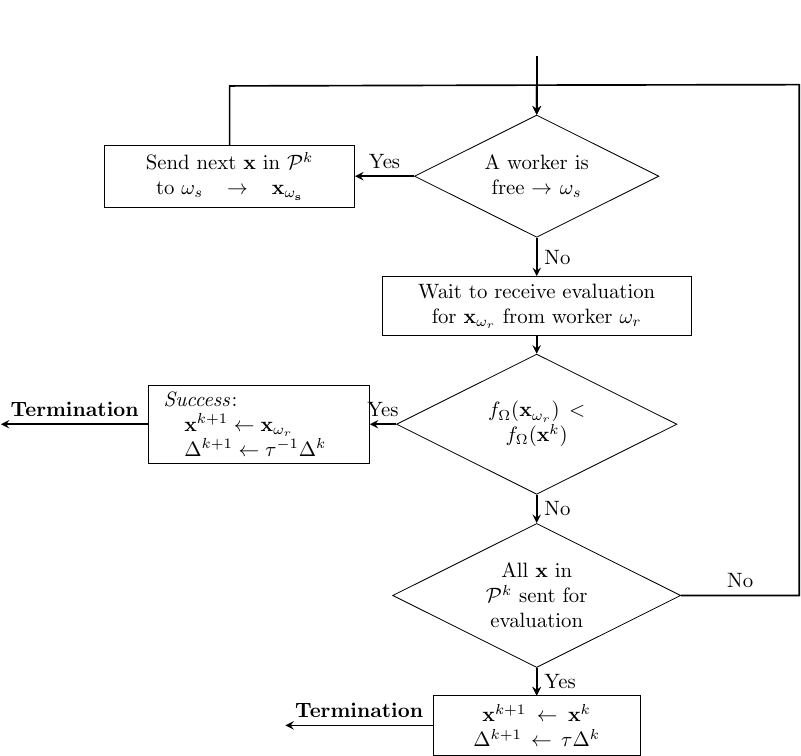}
  \caption{Flow chart of the master process for pMads-A. Details the shaded area in Poll step of \cref{algo:mads}.}
  \label{fig:master_pMadsA}
\end{figure}

As mentioned in~\cite{Le09b}, p\mads-A is a simplified version of the asynchronous parallel pattern search algorithm (\texttt{APPSPACK})~\cite{GrKo06}. The master is permitted to terminate an iteration and to generate new trial points as soon as a new success is obtained by a worker, even if evaluations are still in progress on other workers. 

In the \nomad implementation of p\mads-A, evaluations currently in progress are never cancelled (see \cref{fig:master_pMadsA}).
But once all $\xx$ in $\mathcal{P}^k$ have been sent to workers for evaluation and no success has been obtained among the received evaluations, the current incumbent solution and the frame size are updated. In the event that an evaluation that took place in an iteration prior to the current one finds the best solution, \nomad backtracks and sets this solution as the incumbent one during a cache search step.
Therefore, there is no synchronization barrier between \texttt{MPI} processes during the Poll step. This makes this mode effective for time-varying blackboxes because no \texttt{MPI} processes need to wait for the most time-consuming one to finish evaluating its solution. It is imperative to avoid idling computing resources, as this leads to poor overall performance and induces bottlenecks in a parallel programming context. This assertion is elaborated upon in \cref{sec-tests} where numerical results are provided and discussed to support this claim. The main downside of the asynchronous mode is p\mads-A's non-deterministic behaviour, which limits the reproducibility of the results. We remark that p\mads-A is the default mode used in \nomad when parallelism is enabled. Also, one can note that there still exists a synchronization barrier between the Search and Poll steps that can affect the algorithm effectiveness. 
Notice that, to account for the possibility that evaluations may still be in progress during the \textbf{Termination} step, the stopping criteria of \cref{algo:mads} has to check if all evaluations are terminated.

\subsection{COOP-MADS}
\label{sec-coopmads}
To accelerate the exploration of $\X$ and to prevent a fast convergence of \mads to a local optimum, multiple paths can be generated by running $c$ different instances of \mads in parallel. This can be particularly beneficial for multi-modal problems where the aim is to converge at a global optimum. Each \mads instance is initialized with a unique random seed, which leads to a varied selection of directions from Ortho\mads~\cite{AbAuDeLe09}. Hence, it yields distinct behavioural patterns in trial points generation. 
These \mads instances operate independently and without any synchronization mechanisms between them. To balance the evaluation effort, a dynamic distribution of the global evaluation budget is implemented between instances. Note that a static budget allocation which consists in distributing evaluations equally between instances often results in the resource under-utilization due to the varying convergence rates of individual instances. 

In COOP-\mads, the only information exchange between \mads instances is via the shared cache to prevent redundant evaluations in an instance and between the instances.
The optional cache search could also be enabled to periodically synchronize the incumbent solution between the instances. This option was not used in this work.
Also, for simplicity, each \mads instance runs as a single process to perform the algorithm and the evaluations. The master-worker paradigm used in p\mads is not considered for COOP-\mads. 

\subsection{PSD-MADS}
\label{sec-psdmads}

Solving high-dimensional BBO problems is a computationally demanding task.
PSD-\mads~\cite{AuDeLe07} stands out as an algorithm tailored specifically to large-scale BBO problems. It specifically leverages a parallel space decomposition approach to increase \mads' scalability.
The core concept of PSD-\mads is centred around the manager process generating optimization subproblems $\PP_{p}$ by randomly selecting variable subspaces~$\mathcal{N}_p$ from the original extensive variable space $\mathcal{N}$, i.e., $\mathcal{N}_p \subseteq \mathcal{N}$. These subproblems are then distributed among workers by the manager. Each worker applies \mads on their assigned subproblem and aims to improve the shared incumbent solution $\mathbf{x}^*$. Upon completing its task a worker receives a new subproblem $\PP_p$ from the manager. A subproblem is defined by the incumbent solution given as an initial point, the initial and minimum mesh sizes $\Delta^p_0$ and $\Delta^p_{\min}$ and the fixed variables in $\mathcal{N}$ to define the variable subspace $\mathcal{N}_p$.

Note that the PSD-\mads operates asynchronously as no synchronization step is present. The algorithm introduces four roles, namely, the manager, the pollster, the regular worker, and the cache server. The pollster is a special worker which solves the original problem~\eqref{pb-def} with a single-poll direction \mads algorithm instead of the conventional $2n$ or $n+1$ poll directions. The pollster is crucial to ensure the convergence of PSD-\mads. The manager assumes the responsibility of selecting the optimization subproblems $\PP_p$ and performs adjustments to the main mesh and continuously updates the value of the incumbent solution. The cache server stores all evaluated points to avoid unnecessary, multiple, and expensive function evaluations. Finally, a regular worker is a \mads instance that solves a subproblem $\PP_p$.

    \begin{figure}[ht!]
     \centering
     
        \begin{tikzpicture}[every text node part/.style={font=\tiny}]
        \node[draw, rounded corners=4pt, fill=darkgray!30](master) {\textrm{\begin{tabular}{c}
             Master + \\ 
             Cache Server
        \end{tabular}}};

        \node[draw, rounded corners=4pt, fill=purple!50](w4) at ([shift={(-35mm,0mm)}]master.west){\textrm{\textrm{\begin{tabular}{c}
             \textbf{worker 4}\\ 
             $ns$ variables\\
             $2ns$ directions
        \end{tabular}}}};

        \draw[-stealth]([shift={(-1mm,-1mm)}]master.west)--([shift={(1mm,-1mm)}]w4.east)node[midway,below]{ \tiny subproblem data};
        \draw[stealth-]
        ([shift={(-1mm,1mm)}]master.west)--([shift={(1mm,1mm)}]w4.east)node[midway,above]{ \tiny optimization result};

        \node[draw, rounded corners=4pt, fill=purple!50](w3) at (180 + 50:30mm){\textrm{\textrm{\begin{tabular}{c}
             \textbf{worker 3}\\ 
             $ns$ variables\\
             $2ns$ directions
        \end{tabular}}}};

        \draw[-stealth]([shift={(0mm,1mm)}]w3.north)--([shift={(1mm,-1mm)}]master.210) node[midway,above, sloped]{};
        \draw[stealth-]([shift={(3mm,1mm)}]w3.north)--([shift={(4mm,-1mm)}]master.210) node[midway,above, sloped]{};

        \node[draw, rounded corners=4pt, fill=purple!50](w1) at (50:30mm){\textrm{\textrm{\begin{tabular}{c}
             \textbf{worker 1}\\ 
             $ns$ variables\\
             $2ns$ directions
        \end{tabular}}}};

        \node[draw, rounded corners=4pt, fill=purple!50](w2) at (180 + 130:30mm){\textrm{\textrm{\begin{tabular}{c}
             \textbf{worker 2}\\ 
             $ns$ variables\\
             $2ns$ directions
        \end{tabular}}}};

        \draw[-stealth]([shift={(0mm,1mm)}]w2.north)--([shift={(1mm,-1mm)}]master.320) node[midway,above, sloped]{};
        \draw[stealth-]([shift={(3mm,1mm)}]w2.north)--([shift={(4mm,-1mm)}]master.320) node[midway,above, sloped]{};
        
        \node[draw, rounded corners=4pt, fill=purple!50](poll) at (130:30mm){\textrm{\textrm{\begin{tabular}{c}
             \textbf{pollster} \\ 
             $n$ variables\\
             1 direction
        \end{tabular}}}};
        \draw[-stealth]([shift={(0mm,-1mm)}]poll.south)--([shift={(1mm,1mm)}]master.150) node[midway,above, sloped]{};
        \draw[stealth-]([shift={(3mm,-1mm)}]poll.south)--([shift={(4mm, 1mm)}]master.150) node[midway,above, sloped]{};

        \draw[-stealth]([shift={(0mm,-1mm)}]w1.south)--([shift={(1mm,1mm)}]master.40) node[midway,above, sloped]{};
        \draw[stealth-]([shift={(3mm,-1mm)}]w1.south)--([shift={(4mm, 1mm)}]master.40) node[midway,above, sloped]{};
        
        \end{tikzpicture}
    \caption{Master-worker communication organization for PSD-\mads, where for each worker $p$, the optimization result corresponds to the final mesh size $\Delta^p_\text{end}$ and the final solution $\xx_p$, while the subproblem data corresponds to the starting solution $\xx_0$, the set of variables $\mathcal{N}_p$, and the initial and minimum mesh sizes $\Delta^p_0$ and $\Delta^p_{\min}$.}
        \label{fig:psdmads}
    \end{figure}
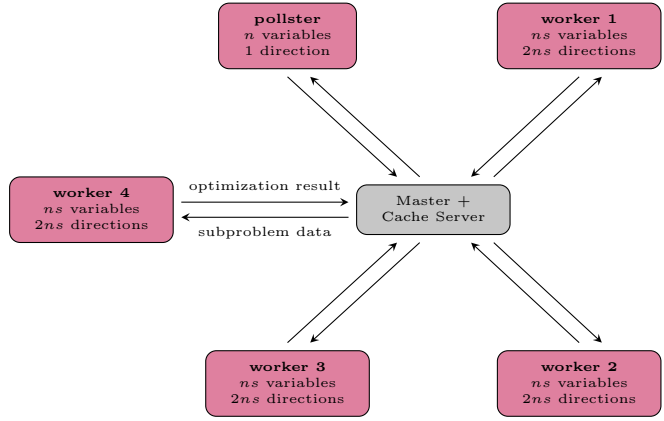

\section{Computational study}
\label{sec-tests}
In this section, we conduct numerical experiments to evaluate the performance of the described algorithms. We use different sets of test problems each of which are targeting specific blackbox characteristics, to make recommendations on when an algorithm should be used. All the experiments are run on Intel Core i7-12700 @ 2.10GHz processors using Version \texttt{3.9.1} of \nomad~\cite{Le09b} with \texttt{MPICH 4.3.0}. We disable the use of quadratic models because this feature is not available in PSD-\mads. %
Our tests use Ortho\mads $2n$ due to the lack of an Ortho\mads $n + 1$ implementation in p\mads-A and PSD-\mads. The number of parallel resources used is set according to the dimension of the given problem and the number of parallel workers that can compute an evaluation given an algorithm. Hence, depending on the algorithm, the number of parallel workers changes. For p\mads, a single process out of the total number provided to \texttt{MPI} is used as the master, leading to the use of $n+1$ processes. For COOP-\mads, we consider $n+1$ \mads instance processes for the algorithm. For PSD-\mads, we use a master and a cache server that do not compute evaluations, leading then to a total of $n+2$ processes.
Blackbox evaluations are performed in batch mode, with input and output files managed by \nomad.

\subsection{Scalability analysis}

The performance of a parallel implementation is mainly characterized by its capacity to reduce a task execution time based on the computing units it uses, which refers to scalability. We use the following two metrics to measure scalability: speedup and efficiency. Let~$S(\eta, p)$ denote the speedup for a task of size $\eta$ using $p$ computing units, defined as:
\begin{equation*}
    S(\eta, p) = \frac{T_1(\eta)}{T_p(\eta)},
\end{equation*}
where $T_1(\eta)$ and $T_p(\eta)$ are the execution times of the sequential and parallel implementations, respectively.
Let the efficiency $E(\eta, p)$ be defined as:
\begin{equation*}
    E(\eta, p) = \frac{S(\eta, p)}{p}.
\end{equation*}
The efficiency relates the speedup to the number of parallel processing units used to achieve the task. 
In our case, the task to parallelize is not predetermined and it is the result of a dynamic aggregation of smaller tasks, namely, the blackbox functions evaluations, as the optimization algorithm proceeds. Hence, efficiency and speedup depend on communication overheads, blackbox evaluation times, and possible synchronization barriers in the algorithm.
Maximum speedup and efficiency occur when communication and synchronization overheads stay small for a relatively large number of computing units $p$.

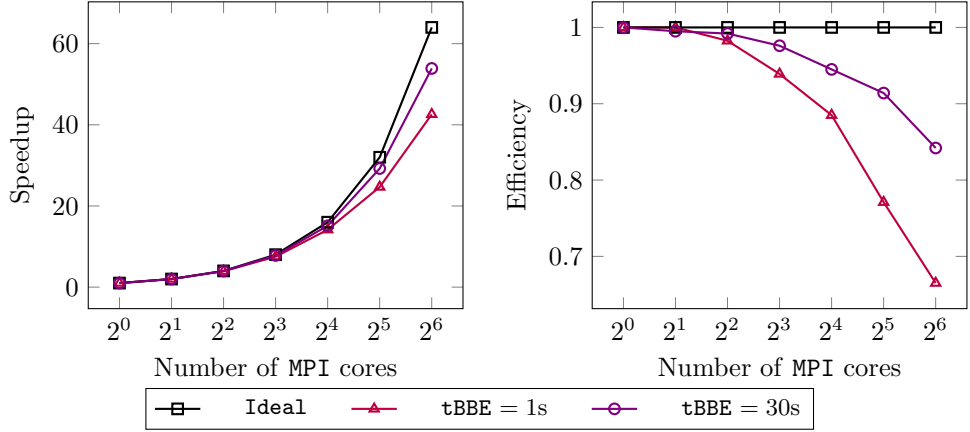
\begin{figure}[ht!]
\begin{subfigure}{0.5\textwidth}
    \centering
    \begin{tikzpicture}
        \begin{axis}[width=\linewidth,xlabel= Number of \texttt{MPI} cores,ylabel= Speedup, xmode=log, log basis x=2, legend cell align={left}, legend pos=north west]
        
        \addplot[purple, mark=triangle, thick]coordinates{ (2^0, 1) (2^1, 2) (2^2, 3.930) (2^3, 7.514) (2^4, 14.153) (2^5, 24.674) (2^6, 42.59)};
        \addplot [black, mark=square, thick] coordinates{ (2^0, 1) (2^1, 2) (2^2, 4) (2^3, 8) (2^4, 16) (2^5, 32) (2^6, 64)};
        \addplot [violet, mark=o, thick] coordinates{ (2^0, 1) (2^1, 1.989) (2^2, 3.966) (2^3, 7.805) (2^4, 15.118) (2^5, 29.236) (2^6, 53.892)};
        
        \end{axis}
    \end{tikzpicture}
\end{subfigure}
\hfill
\begin{subfigure}{0.5\textwidth}
    \centering
    \begin{tikzpicture}
        \begin{axis}[width=\linewidth,xlabel= Number of \texttt{MPI} cores, ylabel= Efficiency, xmode=log, log basis x=2, legend cell align={left}]
        
        \addplot[purple, mark=triangle, thick]coordinates{ (2^0, 1) (2^1, 1) (2^2, 0.9825) (2^3, 0.939) (2^4, 0.885) (2^5, 0.771) (2^6, 0.665)};
        \addplot [black, mark=square, thick] coordinates{ (2^0, 1) (2^1, 1) (2^2, 1) (2^3, 1) (2^4, 1) (2^5, 1) (2^6, 1)};
        \addplot [violet, mark=o, thick] coordinates{ (2^0, 1) (2^1, 0.995) (2^2, 0.992) (2^3, 0.976) (2^4, 0.945) (2^5, 0.914) (2^6, 0.842)};
        
        \end{axis}
    \end{tikzpicture}
\end{subfigure}
\begin{center}
\begin{tikzpicture}
    \begin{axis}[
        hide axis,
        xmin=0, xmax=1,
        ymin=0, ymax=1,
        legend cell align={left},
        legend columns=3,
        legend style={
            at={(0.5, -0.2)},
            anchor=north,
            draw=black,
            font=\small,
            column sep=0.5cm, %
            row sep=4pt, %
            fill=white %
        }
    ]
        \addlegendimage{const plot, black, mark=square, solid, thick}
        \addlegendentry{\texttt{Ideal}}
        
        \addlegendimage{const plot, purple, mark=triangle, thick}
        \addlegendentry{\texttt{tBBE} $=1$s}
        
        \addlegendimage{const plot, violet, mark=o, thick}
        \addlegendentry{\texttt{tBBE} $=30$s}
    \end{axis}
\end{tikzpicture}    
\end{center}
    \caption{Parallel speedup and efficiency of p\mads-A measured over an optimization problem of $32$ variables with a budget of $500$ evaluations and a blackbox evaluation time (t\texttt{BBE}) of $1$s and $30$s. The reported number of \texttt{MPI} cores represents the worker processors.}
    \label{fig:eff_speedup}
\end{figure}

A scalability analysis is performed with a fixed size $\eta$ of 500 blackbox evaluations. A comparison between p\mads-A runs on a blackbox from~\cite{Karmitsa2007}' collection set with different artificial evaluation times (\texttt{tBBE}) is presented in \cref{fig:eff_speedup}. 
These results illustrate that parallel implementation significantly accelerates the computation time. Even at the lowest measured efficiency, a speedup of $43$ is achieved using 64 processors. When subject to an evaluation time of $30$ seconds, the speedup increases from $43$ to $54$.
Hence, for more time-consuming blackbox functions, the parallel performance improves further, as the impact of communication overheads becomes negligible.

\subsection{Heterogeneous test}
\label{sssec:hetero}
Next, we illustrate the impact of the synchronisation barrier over the computation time by comparing p\mads-S and p\mads-A.
The {\sf solar} simulator~\cite{solar_paper} implements in \texttt{C++} a concentrated solar power thermal plant which relies on numerical methods such as Monte Carlo simulation, Newton's method, kernel smoothing, and other iterative methods. The heterogeneous nature of the blackbox comes from combining these methods. The blackbox {\sf solar} offers ten instances, from which we use {\sf solar4.1} because it presents the largest dimension for a single-objective constrained optimization problem, that is, $29$ variables and $16$ constraints. 

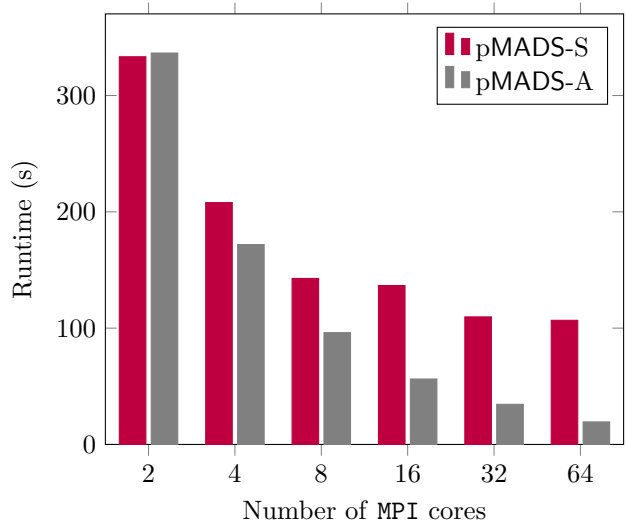
\begin{figure}[ht!]
    \centering
    \begin{tikzpicture}
        \begin{axis}[ybar, ymin=0, xlabel = {Number of \texttt{MPI} cores}, ylabel = {Runtime (s)}, symbolic x coords={2,4,8,16,32,64}, bar width = 10 pt, legend cell align={left}, legend pos=north east, legend entries = {p\mads-S,p\mads-A}]
         \addplot[purple, fill = purple]coordinates {(2, 333.347) (4, 207.959) (8, 142.727) (16, 136.656) (32, 109.678) (64, 106.749)};
        \addplot[gray, fill = gray]coordinates {(2, 336.547) (4, 171.841) (8, 96.152) (16, 56.274) (32, 34.515) (64, 19.388)};
         \end{axis}
    \end{tikzpicture}

\caption{Running time in seconds of p\mads-A and p\mads-S for different numbers of \texttt{MPI} cores when solving {\sf solar4.1} with a $500$ blackbox evaluation budget.
The reported number of \texttt{MPI} cores represents both the master and worker processors.}
\label{fig:solar4_test}
\end{figure}

\cref{fig:solar4_test} illustrates that as more processors are involved, the run time difference between p\mads-S and p\mads-A increases. This is due to the presence of the poll synchronization barrier in p\mads-S because, e.g., it only takes one evaluation with a significantly higher computation time than the others running in parallel to cause substantial computation resource idling. Increasing the number of cores then leads to more evaluations being launched in parallel, thereby increasing the likelihood of encountering a more costly evaluation at each iteration. The potential slow down is particularly sensitive to the blackbox evaluation time distribution but is completely avoided in the case of p\mads-A which steadily decreases its run time.

Due to its clear limitation, for the rest of this work, p\mads-S will be omitted.

\subsection{Benchmarking with data profiles}
Data profiles are used for benchmarking algorithms~\cite{MoWi2009}. Comparisons are done in terms of computation time and number of evaluations scaled with problem dimensions $n_p$. The data profile of an algorithm is a series of values $d_a(k)$ that gives the proportion of $\tau$-solved problems, with $k$ the groups of $n_p + 1$ function evaluations or the time~\cite{G-2025-36}. The accuracy of an evaluation point is obtained at a given time or evaluation number for a given problem using the best and the initial objective values. 

In this work, we have considered the cross-instances best objective value $f^*$ obtained by all algorithms on all run instances of a given problem. The initial points of problems are given and fixed. Different run instances are obtained by providing fixed seeds to run algorithms.

\subsection{Mor\'e-Wild tests}
In this section, the parallel versions of \mads, are compared on the Mor\'e-Wild~\cite{MoWi2009} problems.
The Mor\'e-Wild collection is a well-know test set for benchmarking unconstrained \dfo algorithms in low dimensions. We select the SMOOTH subset which provides $53$ unconstrained optimization problems where the objective function is twice continuously differentiable. For each parallel version of \mads and each problem, runs are conducted with $10$ different random seeds that affect the trial points generation. This gives a total of $530$ run instances per algorithm tested.

\begin{figure}[ht!]
    \centering
    \includegraphics[width=\textwidth]{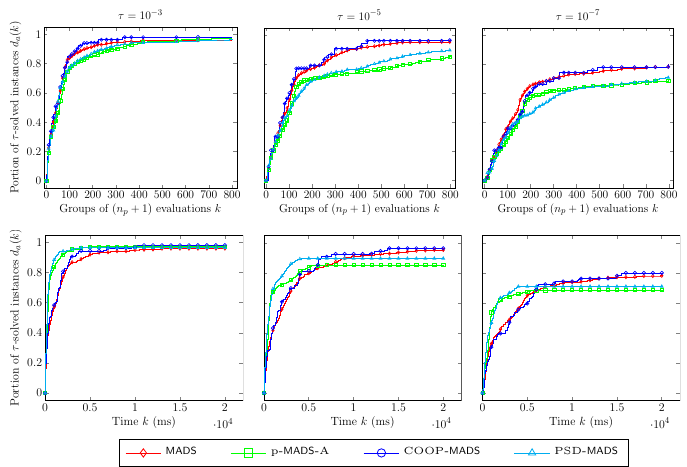}
\caption{Data profiles on the Mor\'e-Wild SMOOTH problem set.
}
\label{fig:morewilde}
 \end{figure}

In~\cref{fig:morewilde}, for $\tau = 10^{-5}$ and below, \mads and COOP-\mads solve more problems than the p\mads and PSD-\mads for the given evaluation budgets. We note that evaluation profiles for p\mads and PSD-\mads are still increasing towards the end. This probably indicates that more problem could have been solved by giving a higher evaluation budget. 
When looking at portion of solved instances with respect to time, p\mads and PSD-\mads achieve their best results faster, with flatter profiles. 

\subsection{Constrained test}
Lastly, we consider a set of continuous constrained problems and use parallel versions of \mads with progressive barrier~\cite{AuDe09a}.
The progressive barrier offers a more sophisticated approach to handling inequality constraints compared to the extreme barrier described in \cref{algo:mads}.
As the default strategy in \nomad, the progressive barrier influences both the detection of success and the management of incumbents, allowing for the identification of two incumbent solutions; a feasible and an infeasible one. Despite its added complexity, it has demonstrated greater effectiveness than the extreme barrier, while maintaining comparable parallelization capabilities.

\begin{figure}[ht!]
    \centering
    \includegraphics[width=\textwidth]{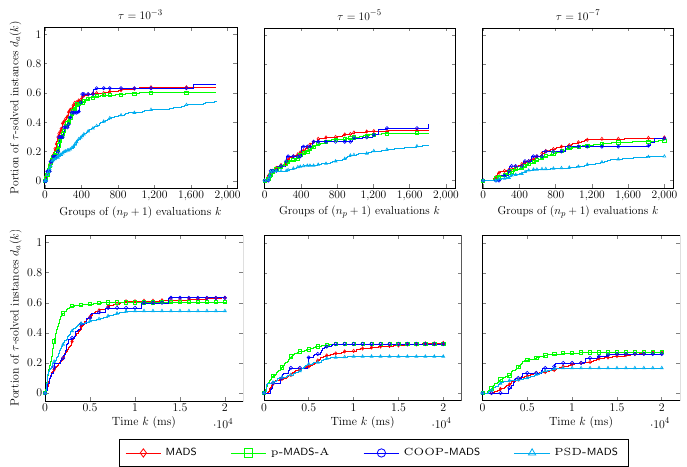}
\caption{Data profiles on the inequality-constrained problem set of \cref{tab:const}.}
\label{fig:constrained}
\end{figure}

\Cref{tab:const} lists and references the problem instances. Ten seeds are again considered to increase the problem collection. Data profiles for different tolerances $\tau$ are compared in \cref{fig:constrained}. 
PSD-\mads do not perform well compared with the other algorithms. The strategy of space decomposition seems not well adapted on small dimension problems with inequality constraints. p\mads-A performance is higher on smaller computation time and as good as COOP-\mads and \mads in terms of evaluations.   

\begin{table*}
\caption{Description of the selected continuous constrained test set. The superscript~$^\S$ denotes the use of multiple starting points for the given problem. }
{\small
\begin{subtable}{0.4\linewidth}
\centering
\begin{tabular}{|c l c@{\hskip.8mm} c@{\hskip.8mm} c@{\hskip.8mm} c|}
\hline
\# & Name & Ref. & $n$ & $m$ & Bnds\\ [0.5ex] 
\hline\hline
1 & CHENWANG-F2$^\S$ &\cite{ChWa2010b}     &  8 & 6 & yes \\
2 & CHENWANG-F3$^\S$ &\cite{ChWa2010b}     & 10 & 8 & yes \\
3 & CRESCENT         &\cite{AuDe09a}       & 10 & 2 & yes \\
4 & DISK             &\cite{AuDe09a}       & 10 & 1 & no \\
5 & FLOUDAS          &{\cite{Floudas1995}} &  3 & 3 & yes \\
6 & G2               &\cite{HeFu06}        & 10 & 2 & yes \\
7 & G9               &\cite{HeFu06}        &  7 & 4 & yes \\
8 & HS83             &\cite{HoSc1981}      &  5 & 6 & yes \\
\hline
\end{tabular}
\end{subtable}
\hspace*{1.25cm}
\begin{subtable}{0.4\linewidth}
\centering
\begin{tabular}{|c l c@{\hskip.8mm} c@{\hskip.8mm} c@{\hskip.8mm} c|}
\hline
\# & Name & Ref. & $n$ & $m$ & Bnds\\ [0.5ex] 
\hline\hline
9 & HS108            &\cite{HoSc1981}    &  9 & 13 & yes \\
10 & HS114           &\cite{HoSc1981}    &  9 &  6 & yes \\
11 & MAD6$^\S$       &\cite{ChWa2010b}   &  5 &  7 & yes \\
12 & OPTENG-RBF      &\cite{AuDe09a}     &  3 &  4 & yes \\
13 & PENTAGON        &\cite{LuVl00}      &  6 & 15 & no \\
14 & SPRING$^\S$     &\cite{RodRenWat98} &  3 &  4 & yes \\
15 & TAOWANG-F2$^\S$ &\cite{TaoWan08}    &  7 &  4 & yes \\
16 & ZHAOWANG-F5     &\cite{ChWa2010b}   & 13 &  9 & yes \\
\hline
\end{tabular}
\end{subtable}
\label{tab:const}
}
\end{table*}

\subsection{Discussion}

The previous subsections numerically illustrates the effectiveness of parallel variants of the \mads algorithm.
The results highlight the critical importance of employing an asynchronous parallel strategy when subject to heterogeneous evaluation times among blackbox functions. Lastly, the p\mads implementations demonstrates strong scalability performance with its efficient utilization of up to $64$ computing cores. This in turn confirms its suitability for parallel computing contexts.

\section{Closing remarks}\label{conclusion}

In this work, we discuss and numerically compare the parallel extensions of the mesh adaptive direct search (\mads) algorithm for blackbox optimization (BBO). First, variants with a synchronization barrier, parallel \mads-Synchronous (p\mads-S) and without, parallel \mads-Asynchronous (p\mads-A), between evaluations are presented. While omitting the barrier improves computational efficiency, it also leads to the non-deterministic behaviour of p\mads-A. Next, the cooperative-\mads (COOP-\mads) algorithm, where several \mads instances are run in parallel with a shared cache, is outlined. Lastly, we survey the parallel space decomposition \mads (PSD-\mads), which specifically targets high-dimensional BBO problems by employing \mads on variable subspaces. The aforementioned extensions are designed such that the local convergence analysis of \mads holds for them as well. Finally, the benefits of each variant, e.g., in terms of scalability, efficiency, running time and/or number of evaluations, is exemplified in numerical simulations. In both unconstrained and constrained problems, p-\mads-A stands out as a time-effective extension for BBO on multiple cores.
In the future, our study is to be extended to the \nomad~4 implementation of \mads. While, the resolution performance should be similar, the computation time and resource scalability may be affected by the its shared memory implementation.

\section*{Acknowledgments}
This work is supported by the NSERC Alliance-Mitacs Accelerate grant  ALLRP~571311-21 (``Optimization of future energy systems'') in collaboration with Hydro-Qu\'ebec.

\section*{Data availability statement}
The open-source solver \nomad~\cite{Le09b} is available at
\href{https://www.gerad.ca/en/software/nomad/}{\tt https://www.gerad.ca/en/software/nomad/}.

\section*{Use of AI statement}
The authors used AI-assisted tools (ChatGPT, OpenAI) solely to improve the language and readability of the manuscript; all scientific content and conclusions are entirely the authors’ own.

\bibliography{bibliography}
\pdfbookmark[1]{References}{sec-refs}
\label{sec-refs}

\end{document}